\documentclass[11pt,oneside,article]{memoir}

\settrims{0pt}{0pt} %
\settypeblocksize{*}{35.5pc}{*} %
\setlrmargins{*}{*}{1} %
\setulmarginsandblock{1in}{1in}{*} %
\setheadfoot{\onelineskip}{1.5\onelineskip} %
\setheaderspaces{*}{1.5\onelineskip}{*} %
\checkandfixthelayout

\usepackage{mathtools}
\usepackage{amsthm}
\usepackage{amssymb}
\usepackage[T1]{fontenc}
\usepackage{newpxtext}
\usepackage[varg,bigdelims]{newpxmath}
\usepackage{stmaryrd}
\usepackage[cal=euler,scr=rsfso]{mathalfa}
\usepackage{dutchcal}
\usepackage[draft]{minted}
\usepackage{bm}
\usepackage{inputenc}
\usepackage{microtype}
\usepackage[usenames,dvipsnames]{xcolor}
\usepackage{enumitem}
\usepackage{xcolor}
\usepackage{booktabs}
\usepackage{tikz}
\usepackage[prefix=tikzsym]{tikzsymbols}
\usepackage{todonotes}
\RequirePackage{hyperref}
\hypersetup{colorlinks,linkcolor={blue!50!black},citecolor={red!50!black},urlcolor={blue!80!black}}

\usepackage[backend=biber, maxbibnames = 10, style = alphabetic]{biblatex}
  \newtheorem{theorem}{Theorem}[chapter]
  \newtheorem{proposition}[theorem]{Proposition}

  \theoremstyle{definition}
  \newtheorem{definition}[theorem]{Definition}
  \newtheorem{notation}[theorem]{Notation}
  
  \newtheorem*{axiom*}{Axiom}
  
  \theoremstyle{remark}
  \newtheorem{example}[theorem]{Example}
  \newtheorem{remark}[theorem]{Remark}

\makeatletter
\def\PYG@reset{\let\PYG@it=\relax \let\PYG@bf=\relax%
    \let\PYG@ul=\relax \let\PYG@tc=\relax%
    \let\PYG@bc=\relax \let\PYG@ff=\relax}
\def\PYG@tok#1{\csname PYG@tok@#1\endcsname}
\def\PYG@toks#1+{\ifx\relax#1\empty\else%
    \PYG@tok{#1}\expandafter\PYG@toks\fi}
\def\PYG@do#1{\PYG@bc{\PYG@tc{\PYG@ul{%
    \PYG@it{\PYG@bf{\PYG@ff{#1}}}}}}}
\def\PYG#1#2{\PYG@reset\PYG@toks#1+\relax+\PYG@do{#2}}

\expandafter\def\csname PYG@tok@gd\endcsname{\def\PYG@tc##1{\textcolor[rgb]{0.63,0.00,0.00}{##1}}}
\expandafter\def\csname PYG@tok@gu\endcsname{\let\PYG@bf=\textbf\def\PYG@tc##1{\textcolor[rgb]{0.50,0.00,0.50}{##1}}}
\expandafter\def\csname PYG@tok@gt\endcsname{\def\PYG@tc##1{\textcolor[rgb]{0.00,0.27,0.87}{##1}}}
\expandafter\def\csname PYG@tok@gs\endcsname{\let\PYG@bf=\textbf}
\expandafter\def\csname PYG@tok@gr\endcsname{\def\PYG@tc##1{\textcolor[rgb]{1.00,0.00,0.00}{##1}}}
\expandafter\def\csname PYG@tok@cm\endcsname{\let\PYG@it=\textit\def\PYG@tc##1{\textcolor[rgb]{0.25,0.50,0.50}{##1}}}
\expandafter\def\csname PYG@tok@vg\endcsname{\def\PYG@tc##1{\textcolor[rgb]{0.10,0.09,0.49}{##1}}}
\expandafter\def\csname PYG@tok@vi\endcsname{\def\PYG@tc##1{\textcolor[rgb]{0.10,0.09,0.49}{##1}}}
\expandafter\def\csname PYG@tok@vm\endcsname{\def\PYG@tc##1{\textcolor[rgb]{0.10,0.09,0.49}{##1}}}
\expandafter\def\csname PYG@tok@mh\endcsname{\def\PYG@tc##1{\textcolor[rgb]{0.40,0.40,0.40}{##1}}}
\expandafter\def\csname PYG@tok@cs\endcsname{\let\PYG@it=\textit\def\PYG@tc##1{\textcolor[rgb]{0.25,0.50,0.50}{##1}}}
\expandafter\def\csname PYG@tok@ge\endcsname{\let\PYG@it=\textit}
\expandafter\def\csname PYG@tok@vc\endcsname{\def\PYG@tc##1{\textcolor[rgb]{0.10,0.09,0.49}{##1}}}
\expandafter\def\csname PYG@tok@il\endcsname{\def\PYG@tc##1{\textcolor[rgb]{0.40,0.40,0.40}{##1}}}
\expandafter\def\csname PYG@tok@go\endcsname{\def\PYG@tc##1{\textcolor[rgb]{0.53,0.53,0.53}{##1}}}
\expandafter\def\csname PYG@tok@cp\endcsname{\def\PYG@tc##1{\textcolor[rgb]{0.74,0.48,0.00}{##1}}}
\expandafter\def\csname PYG@tok@gi\endcsname{\def\PYG@tc##1{\textcolor[rgb]{0.00,0.63,0.00}{##1}}}
\expandafter\def\csname PYG@tok@gh\endcsname{\let\PYG@bf=\textbf\def\PYG@tc##1{\textcolor[rgb]{0.00,0.00,0.50}{##1}}}
\expandafter\def\csname PYG@tok@ni\endcsname{\let\PYG@bf=\textbf\def\PYG@tc##1{\textcolor[rgb]{0.60,0.60,0.60}{##1}}}
\expandafter\def\csname PYG@tok@nl\endcsname{\def\PYG@tc##1{\textcolor[rgb]{0.63,0.63,0.00}{##1}}}
\expandafter\def\csname PYG@tok@nn\endcsname{\let\PYG@bf=\textbf\def\PYG@tc##1{\textcolor[rgb]{0.00,0.00,1.00}{##1}}}
\expandafter\def\csname PYG@tok@no\endcsname{\def\PYG@tc##1{\textcolor[rgb]{0.53,0.00,0.00}{##1}}}
\expandafter\def\csname PYG@tok@na\endcsname{\def\PYG@tc##1{\textcolor[rgb]{0.49,0.56,0.16}{##1}}}
\expandafter\def\csname PYG@tok@nb\endcsname{\def\PYG@tc##1{\textcolor[rgb]{0.00,0.50,0.00}{##1}}}
\expandafter\def\csname PYG@tok@nc\endcsname{\let\PYG@bf=\textbf\def\PYG@tc##1{\textcolor[rgb]{0.00,0.00,1.00}{##1}}}
\expandafter\def\csname PYG@tok@nd\endcsname{\def\PYG@tc##1{\textcolor[rgb]{0.67,0.13,1.00}{##1}}}
\expandafter\def\csname PYG@tok@ne\endcsname{\let\PYG@bf=\textbf\def\PYG@tc##1{\textcolor[rgb]{0.82,0.25,0.23}{##1}}}
\expandafter\def\csname PYG@tok@nf\endcsname{\def\PYG@tc##1{\textcolor[rgb]{0.00,0.00,1.00}{##1}}}
\expandafter\def\csname PYG@tok@si\endcsname{\let\PYG@bf=\textbf\def\PYG@tc##1{\textcolor[rgb]{0.73,0.40,0.53}{##1}}}
\expandafter\def\csname PYG@tok@s2\endcsname{\def\PYG@tc##1{\textcolor[rgb]{0.73,0.13,0.13}{##1}}}
\expandafter\def\csname PYG@tok@nt\endcsname{\let\PYG@bf=\textbf\def\PYG@tc##1{\textcolor[rgb]{0.00,0.50,0.00}{##1}}}
\expandafter\def\csname PYG@tok@nv\endcsname{\def\PYG@tc##1{\textcolor[rgb]{0.10,0.09,0.49}{##1}}}
\expandafter\def\csname PYG@tok@s1\endcsname{\def\PYG@tc##1{\textcolor[rgb]{0.73,0.13,0.13}{##1}}}
\expandafter\def\csname PYG@tok@dl\endcsname{\def\PYG@tc##1{\textcolor[rgb]{0.73,0.13,0.13}{##1}}}
\expandafter\def\csname PYG@tok@ch\endcsname{\let\PYG@it=\textit\def\PYG@tc##1{\textcolor[rgb]{0.25,0.50,0.50}{##1}}}
\expandafter\def\csname PYG@tok@m\endcsname{\def\PYG@tc##1{\textcolor[rgb]{0.40,0.40,0.40}{##1}}}
\expandafter\def\csname PYG@tok@gp\endcsname{\let\PYG@bf=\textbf\def\PYG@tc##1{\textcolor[rgb]{0.00,0.00,0.50}{##1}}}
\expandafter\def\csname PYG@tok@sh\endcsname{\def\PYG@tc##1{\textcolor[rgb]{0.73,0.13,0.13}{##1}}}
\expandafter\def\csname PYG@tok@ow\endcsname{\let\PYG@bf=\textbf\def\PYG@tc##1{\textcolor[rgb]{0.67,0.13,1.00}{##1}}}
\expandafter\def\csname PYG@tok@sx\endcsname{\def\PYG@tc##1{\textcolor[rgb]{0.00,0.50,0.00}{##1}}}
\expandafter\def\csname PYG@tok@bp\endcsname{\def\PYG@tc##1{\textcolor[rgb]{0.00,0.50,0.00}{##1}}}
\expandafter\def\csname PYG@tok@c1\endcsname{\let\PYG@it=\textit\def\PYG@tc##1{\textcolor[rgb]{0.25,0.50,0.50}{##1}}}
\expandafter\def\csname PYG@tok@fm\endcsname{\def\PYG@tc##1{\textcolor[rgb]{0.00,0.00,1.00}{##1}}}
\expandafter\def\csname PYG@tok@o\endcsname{\def\PYG@tc##1{\textcolor[rgb]{0.40,0.40,0.40}{##1}}}
\expandafter\def\csname PYG@tok@kc\endcsname{\let\PYG@bf=\textbf\def\PYG@tc##1{\textcolor[rgb]{0.00,0.50,0.00}{##1}}}
\expandafter\def\csname PYG@tok@c\endcsname{\let\PYG@it=\textit\def\PYG@tc##1{\textcolor[rgb]{0.25,0.50,0.50}{##1}}}
\expandafter\def\csname PYG@tok@mf\endcsname{\def\PYG@tc##1{\textcolor[rgb]{0.40,0.40,0.40}{##1}}}
\expandafter\def\csname PYG@tok@err\endcsname{\def\PYG@bc##1{\setlength{\fboxsep}{0pt}\fcolorbox[rgb]{1.00,0.00,0.00}{1,1,1}{\strut ##1}}}
\expandafter\def\csname PYG@tok@mb\endcsname{\def\PYG@tc##1{\textcolor[rgb]{0.40,0.40,0.40}{##1}}}
\expandafter\def\csname PYG@tok@ss\endcsname{\def\PYG@tc##1{\textcolor[rgb]{0.10,0.09,0.49}{##1}}}
\expandafter\def\csname PYG@tok@sr\endcsname{\def\PYG@tc##1{\textcolor[rgb]{0.73,0.40,0.53}{##1}}}
\expandafter\def\csname PYG@tok@mo\endcsname{\def\PYG@tc##1{\textcolor[rgb]{0.40,0.40,0.40}{##1}}}
\expandafter\def\csname PYG@tok@kd\endcsname{\let\PYG@bf=\textbf\def\PYG@tc##1{\textcolor[rgb]{0.00,0.50,0.00}{##1}}}
\expandafter\def\csname PYG@tok@mi\endcsname{\def\PYG@tc##1{\textcolor[rgb]{0.40,0.40,0.40}{##1}}}
\expandafter\def\csname PYG@tok@kn\endcsname{\let\PYG@bf=\textbf\def\PYG@tc##1{\textcolor[rgb]{0.00,0.50,0.00}{##1}}}
\expandafter\def\csname PYG@tok@cpf\endcsname{\let\PYG@it=\textit\def\PYG@tc##1{\textcolor[rgb]{0.25,0.50,0.50}{##1}}}
\expandafter\def\csname PYG@tok@kr\endcsname{\let\PYG@bf=\textbf\def\PYG@tc##1{\textcolor[rgb]{0.00,0.50,0.00}{##1}}}
\expandafter\def\csname PYG@tok@s\endcsname{\def\PYG@tc##1{\textcolor[rgb]{0.73,0.13,0.13}{##1}}}
\expandafter\def\csname PYG@tok@kp\endcsname{\def\PYG@tc##1{\textcolor[rgb]{0.00,0.50,0.00}{##1}}}
\expandafter\def\csname PYG@tok@w\endcsname{\def\PYG@tc##1{\textcolor[rgb]{0.73,0.73,0.73}{##1}}}
\expandafter\def\csname PYG@tok@kt\endcsname{\def\PYG@tc##1{\textcolor[rgb]{0.69,0.00,0.25}{##1}}}
\expandafter\def\csname PYG@tok@sc\endcsname{\def\PYG@tc##1{\textcolor[rgb]{0.73,0.13,0.13}{##1}}}
\expandafter\def\csname PYG@tok@sb\endcsname{\def\PYG@tc##1{\textcolor[rgb]{0.73,0.13,0.13}{##1}}}
\expandafter\def\csname PYG@tok@sa\endcsname{\def\PYG@tc##1{\textcolor[rgb]{0.73,0.13,0.13}{##1}}}
\expandafter\def\csname PYG@tok@k\endcsname{\let\PYG@bf=\textbf\def\PYG@tc##1{\textcolor[rgb]{0.00,0.50,0.00}{##1}}}
\expandafter\def\csname PYG@tok@se\endcsname{\let\PYG@bf=\textbf\def\PYG@tc##1{\textcolor[rgb]{0.73,0.40,0.13}{##1}}}
\expandafter\def\csname PYG@tok@sd\endcsname{\let\PYG@it=\textit\def\PYG@tc##1{\textcolor[rgb]{0.73,0.13,0.13}{##1}}}

\makeatother

\makeatletter
\def\PYGdefault@reset{\let\PYGdefault@it=\relax \let\PYGdefault@bf=\relax%
    \let\PYGdefault@ul=\relax \let\PYGdefault@tc=\relax%
    \let\PYGdefault@bc=\relax \let\PYGdefault@ff=\relax}
\def\PYGdefault@tok#1{\csname PYGdefault@tok@#1\endcsname}
\def\PYGdefault@toks#1+{\ifx\relax#1\empty\else%
    \PYGdefault@tok{#1}\expandafter\PYGdefault@toks\fi}
\def\PYGdefault@do#1{\PYGdefault@bc{\PYGdefault@tc{\PYGdefault@ul{%
    \PYGdefault@it{\PYGdefault@bf{\PYGdefault@ff{#1}}}}}}}
\def\PYGdefault#1#2{\PYGdefault@reset\PYGdefault@toks#1+\relax+\PYGdefault@do{#2}}

\expandafter\def\csname PYGdefault@tok@gd\endcsname{\def\PYGdefault@tc##1{\textcolor[rgb]{0.63,0.00,0.00}{##1}}}
\expandafter\def\csname PYGdefault@tok@gu\endcsname{\let\PYGdefault@bf=\textbf\def\PYGdefault@tc##1{\textcolor[rgb]{0.50,0.00,0.50}{##1}}}
\expandafter\def\csname PYGdefault@tok@gt\endcsname{\def\PYGdefault@tc##1{\textcolor[rgb]{0.00,0.27,0.87}{##1}}}
\expandafter\def\csname PYGdefault@tok@gs\endcsname{\let\PYGdefault@bf=\textbf}
\expandafter\def\csname PYGdefault@tok@gr\endcsname{\def\PYGdefault@tc##1{\textcolor[rgb]{1.00,0.00,0.00}{##1}}}
\expandafter\def\csname PYGdefault@tok@cm\endcsname{\let\PYGdefault@it=\textit\def\PYGdefault@tc##1{\textcolor[rgb]{0.25,0.50,0.50}{##1}}}
\expandafter\def\csname PYGdefault@tok@vg\endcsname{\def\PYGdefault@tc##1{\textcolor[rgb]{0.10,0.09,0.49}{##1}}}
\expandafter\def\csname PYGdefault@tok@vi\endcsname{\def\PYGdefault@tc##1{\textcolor[rgb]{0.10,0.09,0.49}{##1}}}
\expandafter\def\csname PYGdefault@tok@vm\endcsname{\def\PYGdefault@tc##1{\textcolor[rgb]{0.10,0.09,0.49}{##1}}}
\expandafter\def\csname PYGdefault@tok@mh\endcsname{\def\PYGdefault@tc##1{\textcolor[rgb]{0.40,0.40,0.40}{##1}}}
\expandafter\def\csname PYGdefault@tok@cs\endcsname{\let\PYGdefault@it=\textit\def\PYGdefault@tc##1{\textcolor[rgb]{0.25,0.50,0.50}{##1}}}
\expandafter\def\csname PYGdefault@tok@ge\endcsname{\let\PYGdefault@it=\textit}
\expandafter\def\csname PYGdefault@tok@vc\endcsname{\def\PYGdefault@tc##1{\textcolor[rgb]{0.10,0.09,0.49}{##1}}}
\expandafter\def\csname PYGdefault@tok@il\endcsname{\def\PYGdefault@tc##1{\textcolor[rgb]{0.40,0.40,0.40}{##1}}}
\expandafter\def\csname PYGdefault@tok@go\endcsname{\def\PYGdefault@tc##1{\textcolor[rgb]{0.53,0.53,0.53}{##1}}}
\expandafter\def\csname PYGdefault@tok@cp\endcsname{\def\PYGdefault@tc##1{\textcolor[rgb]{0.74,0.48,0.00}{##1}}}
\expandafter\def\csname PYGdefault@tok@gi\endcsname{\def\PYGdefault@tc##1{\textcolor[rgb]{0.00,0.63,0.00}{##1}}}
\expandafter\def\csname PYGdefault@tok@gh\endcsname{\let\PYGdefault@bf=\textbf\def\PYGdefault@tc##1{\textcolor[rgb]{0.00,0.00,0.50}{##1}}}
\expandafter\def\csname PYGdefault@tok@ni\endcsname{\let\PYGdefault@bf=\textbf\def\PYGdefault@tc##1{\textcolor[rgb]{0.60,0.60,0.60}{##1}}}
\expandafter\def\csname PYGdefault@tok@nl\endcsname{\def\PYGdefault@tc##1{\textcolor[rgb]{0.63,0.63,0.00}{##1}}}
\expandafter\def\csname PYGdefault@tok@nn\endcsname{\let\PYGdefault@bf=\textbf\def\PYGdefault@tc##1{\textcolor[rgb]{0.00,0.00,1.00}{##1}}}
\expandafter\def\csname PYGdefault@tok@no\endcsname{\def\PYGdefault@tc##1{\textcolor[rgb]{0.53,0.00,0.00}{##1}}}
\expandafter\def\csname PYGdefault@tok@na\endcsname{\def\PYGdefault@tc##1{\textcolor[rgb]{0.49,0.56,0.16}{##1}}}
\expandafter\def\csname PYGdefault@tok@nb\endcsname{\def\PYGdefault@tc##1{\textcolor[rgb]{0.00,0.50,0.00}{##1}}}
\expandafter\def\csname PYGdefault@tok@nc\endcsname{\let\PYGdefault@bf=\textbf\def\PYGdefault@tc##1{\textcolor[rgb]{0.00,0.00,1.00}{##1}}}
\expandafter\def\csname PYGdefault@tok@nd\endcsname{\def\PYGdefault@tc##1{\textcolor[rgb]{0.67,0.13,1.00}{##1}}}
\expandafter\def\csname PYGdefault@tok@ne\endcsname{\let\PYGdefault@bf=\textbf\def\PYGdefault@tc##1{\textcolor[rgb]{0.82,0.25,0.23}{##1}}}
\expandafter\def\csname PYGdefault@tok@nf\endcsname{\def\PYGdefault@tc##1{\textcolor[rgb]{0.00,0.00,1.00}{##1}}}
\expandafter\def\csname PYGdefault@tok@si\endcsname{\let\PYGdefault@bf=\textbf\def\PYGdefault@tc##1{\textcolor[rgb]{0.73,0.40,0.53}{##1}}}
\expandafter\def\csname PYGdefault@tok@s2\endcsname{\def\PYGdefault@tc##1{\textcolor[rgb]{0.73,0.13,0.13}{##1}}}
\expandafter\def\csname PYGdefault@tok@nt\endcsname{\let\PYGdefault@bf=\textbf\def\PYGdefault@tc##1{\textcolor[rgb]{0.00,0.50,0.00}{##1}}}
\expandafter\def\csname PYGdefault@tok@nv\endcsname{\def\PYGdefault@tc##1{\textcolor[rgb]{0.10,0.09,0.49}{##1}}}
\expandafter\def\csname PYGdefault@tok@s1\endcsname{\def\PYGdefault@tc##1{\textcolor[rgb]{0.73,0.13,0.13}{##1}}}
\expandafter\def\csname PYGdefault@tok@dl\endcsname{\def\PYGdefault@tc##1{\textcolor[rgb]{0.73,0.13,0.13}{##1}}}
\expandafter\def\csname PYGdefault@tok@ch\endcsname{\let\PYGdefault@it=\textit\def\PYGdefault@tc##1{\textcolor[rgb]{0.25,0.50,0.50}{##1}}}
\expandafter\def\csname PYGdefault@tok@m\endcsname{\def\PYGdefault@tc##1{\textcolor[rgb]{0.40,0.40,0.40}{##1}}}
\expandafter\def\csname PYGdefault@tok@gp\endcsname{\let\PYGdefault@bf=\textbf\def\PYGdefault@tc##1{\textcolor[rgb]{0.00,0.00,0.50}{##1}}}
\expandafter\def\csname PYGdefault@tok@sh\endcsname{\def\PYGdefault@tc##1{\textcolor[rgb]{0.73,0.13,0.13}{##1}}}
\expandafter\def\csname PYGdefault@tok@ow\endcsname{\let\PYGdefault@bf=\textbf\def\PYGdefault@tc##1{\textcolor[rgb]{0.67,0.13,1.00}{##1}}}
\expandafter\def\csname PYGdefault@tok@sx\endcsname{\def\PYGdefault@tc##1{\textcolor[rgb]{0.00,0.50,0.00}{##1}}}
\expandafter\def\csname PYGdefault@tok@bp\endcsname{\def\PYGdefault@tc##1{\textcolor[rgb]{0.00,0.50,0.00}{##1}}}
\expandafter\def\csname PYGdefault@tok@c1\endcsname{\let\PYGdefault@it=\textit\def\PYGdefault@tc##1{\textcolor[rgb]{0.25,0.50,0.50}{##1}}}
\expandafter\def\csname PYGdefault@tok@fm\endcsname{\def\PYGdefault@tc##1{\textcolor[rgb]{0.00,0.00,1.00}{##1}}}
\expandafter\def\csname PYGdefault@tok@o\endcsname{\def\PYGdefault@tc##1{\textcolor[rgb]{0.40,0.40,0.40}{##1}}}
\expandafter\def\csname PYGdefault@tok@kc\endcsname{\let\PYGdefault@bf=\textbf\def\PYGdefault@tc##1{\textcolor[rgb]{0.00,0.50,0.00}{##1}}}
\expandafter\def\csname PYGdefault@tok@c\endcsname{\let\PYGdefault@it=\textit\def\PYGdefault@tc##1{\textcolor[rgb]{0.25,0.50,0.50}{##1}}}
\expandafter\def\csname PYGdefault@tok@mf\endcsname{\def\PYGdefault@tc##1{\textcolor[rgb]{0.40,0.40,0.40}{##1}}}
\expandafter\def\csname PYGdefault@tok@err\endcsname{\def\PYGdefault@bc##1{\setlength{\fboxsep}{0pt}\fcolorbox[rgb]{1.00,0.00,0.00}{1,1,1}{\strut ##1}}}
\expandafter\def\csname PYGdefault@tok@mb\endcsname{\def\PYGdefault@tc##1{\textcolor[rgb]{0.40,0.40,0.40}{##1}}}
\expandafter\def\csname PYGdefault@tok@ss\endcsname{\def\PYGdefault@tc##1{\textcolor[rgb]{0.10,0.09,0.49}{##1}}}
\expandafter\def\csname PYGdefault@tok@sr\endcsname{\def\PYGdefault@tc##1{\textcolor[rgb]{0.73,0.40,0.53}{##1}}}
\expandafter\def\csname PYGdefault@tok@mo\endcsname{\def\PYGdefault@tc##1{\textcolor[rgb]{0.40,0.40,0.40}{##1}}}
\expandafter\def\csname PYGdefault@tok@kd\endcsname{\let\PYGdefault@bf=\textbf\def\PYGdefault@tc##1{\textcolor[rgb]{0.00,0.50,0.00}{##1}}}
\expandafter\def\csname PYGdefault@tok@mi\endcsname{\def\PYGdefault@tc##1{\textcolor[rgb]{0.40,0.40,0.40}{##1}}}
\expandafter\def\csname PYGdefault@tok@kn\endcsname{\let\PYGdefault@bf=\textbf\def\PYGdefault@tc##1{\textcolor[rgb]{0.00,0.50,0.00}{##1}}}
\expandafter\def\csname PYGdefault@tok@cpf\endcsname{\let\PYGdefault@it=\textit\def\PYGdefault@tc##1{\textcolor[rgb]{0.25,0.50,0.50}{##1}}}
\expandafter\def\csname PYGdefault@tok@kr\endcsname{\let\PYGdefault@bf=\textbf\def\PYGdefault@tc##1{\textcolor[rgb]{0.00,0.50,0.00}{##1}}}
\expandafter\def\csname PYGdefault@tok@s\endcsname{\def\PYGdefault@tc##1{\textcolor[rgb]{0.73,0.13,0.13}{##1}}}
\expandafter\def\csname PYGdefault@tok@kp\endcsname{\def\PYGdefault@tc##1{\textcolor[rgb]{0.00,0.50,0.00}{##1}}}
\expandafter\def\csname PYGdefault@tok@w\endcsname{\def\PYGdefault@tc##1{\textcolor[rgb]{0.73,0.73,0.73}{##1}}}
\expandafter\def\csname PYGdefault@tok@kt\endcsname{\def\PYGdefault@tc##1{\textcolor[rgb]{0.69,0.00,0.25}{##1}}}
\expandafter\def\csname PYGdefault@tok@sc\endcsname{\def\PYGdefault@tc##1{\textcolor[rgb]{0.73,0.13,0.13}{##1}}}
\expandafter\def\csname PYGdefault@tok@sb\endcsname{\def\PYGdefault@tc##1{\textcolor[rgb]{0.73,0.13,0.13}{##1}}}
\expandafter\def\csname PYGdefault@tok@sa\endcsname{\def\PYGdefault@tc##1{\textcolor[rgb]{0.73,0.13,0.13}{##1}}}
\expandafter\def\csname PYGdefault@tok@k\endcsname{\let\PYGdefault@bf=\textbf\def\PYGdefault@tc##1{\textcolor[rgb]{0.00,0.50,0.00}{##1}}}
\expandafter\def\csname PYGdefault@tok@se\endcsname{\let\PYGdefault@bf=\textbf\def\PYGdefault@tc##1{\textcolor[rgb]{0.73,0.40,0.13}{##1}}}
\expandafter\def\csname PYGdefault@tok@sd\endcsname{\let\PYGdefault@it=\textit\def\PYGdefault@tc##1{\textcolor[rgb]{0.73,0.13,0.13}{##1}}}

\makeatother

\usetikzlibrary{
	cd,
	math,
	decorations.markings,
	decorations.pathreplacing,
	decorations.pathmorphing,
	positioning,
	arrows.meta,
	circuits.logic.US,
	shapes,
	calc,
	fit,
	quotes}
	
\tikzset{trees/.style={
	inner sep=0, 
	minimum width=0, 
	minimum height=0,
	level distance=.5cm, 
	sibling distance=.5cm,
	edge from parent/.style={shorten <= -4pt, draw, ->}
	}
}

  \tikzset{
     oriented WD/.style={%
        every to/.style={out=0,in=180,draw},
        label/.style={
           font=\everymath\expandafter{\the\everymath\scriptstyle},
           inner sep=0pt,
           node distance=2pt and -2pt},
        semithick,
        node distance=1 and 1,
        decoration={markings, mark=at position \stringdecpos with \stringdec},
        ar/.style={postaction={decorate}},
        execute at begin picture={\tikzset{
           x=\bbx, y=\bby,
           }}
        },
     string decoration/.store in=\stringdec,
     string decoration={\arrow{stealth};},
     string decoration pos/.store in=\stringdecpos,
     string decoration pos=.7,
     bbx/.store in=\bbx,
     bbx = 1.5cm,
     bby/.store in=\bby,
     bby = 1.5ex,
     bb port sep/.store in=\bbportsep,
     bb port sep=1.5,
     bb port length/.store in=\bbportlen,
     bb port length=4pt,
     bb penetrate/.store in=\bbpenetrate,
     bb penetrate=0,
     bb min width/.store in=\bbminwidth,
     bb min width=1cm,
     bb rounded corners/.store in=\bbcorners,
     bb rounded corners=2pt,
     bb small/.style={bb port sep=1, bb port length=2.5pt, bbx=.4cm, bb min width=.4cm, 
bby=.7ex},
		 bb medium/.style={bb port sep=1, bb port length=2.5pt, bbx=.4cm, bb min width=.4cm, 
bby=.9ex},
     bb/.code 2 args={%
        \pgfmathsetlengthmacro{\bbheight}{\bbportsep * (max(#1,#2)+1) * \bby}
        \pgfkeysalso{draw,minimum height=\bbheight,minimum width=\bbminwidth,outer 
sep=0pt,
           rounded corners=\bbcorners,thick,
           prefix after command={\pgfextra{\let\fixname\tikzlastnode}},
           append after command={\pgfextra{\draw
              \ifnum #1=0{} \else foreach \i in {1,...,#1} {
                 ($(\fixname.north west)!{\i/(#1+1)}!(\fixname.south west)$) +(-
\bbportlen,0) 
  coordinate (\fixname_in\i) -- +(\bbpenetrate,0) coordinate (\fixname_in\i')}\fi 
              \ifnum #2=0{} \else foreach \i in {1,...,#2} {
                 ($(\fixname.north east)!{\i/(#2+1)}!(\fixname.south east)$) +(-
\bbpenetrate,0) 
  coordinate (\fixname_out\i') -- +(\bbportlen,0) coordinate (\fixname_out\i)}\fi;
           }}}
     },
     bb name/.style={append after command={\pgfextra{\node[anchor=north] at 
(\fixname.north) {#1};}}}
  }

\DeclareMathOperator{\cod}{cod}
\DeclareMathOperator*{\colim}{colim}

\DeclareMathOperator{\ob}{Ob}

\newcommand{\Cat}[1]{\mathsf{#1}}

\newcommand{\fun}[1]{\mathrm{#1}}

\newcommand{\tn}{\textnormal}

\newcommand{\rr}{\mathbb{R}}
\newcommand{\smset}{\Cat{Set}}
\newcommand{\smcat}{\Cat{Cat}}

\newcommand{\poly}{\Cat{Poly}}
\newcommand{\comon}{\Cat{Comon}}

\newcommand{\yon}{\mathcal{y}}

\newcommand{\lmo}[2][over]{\ifthenelse{\equal{#1}{over}}{\overset{#2}{\bullet}}{\underset{#2}{\bullet}}}

\newcommand{\circpow}[1]{^{\circ #1}}

\newcommand{\id}{\mathrm{id}}

\newcommand{\To}[1]{\xrightarrow{#1}}

\newcommand{\from}{\leftarrow}
\newcommand{\From}[1]{\xleftarrow{#1}}

\newcommand{\op}{^{\tn{op}}}

\newcommand{\mdot}{{\cdot}}
\newcommand{\inv}{^{-1}}

\newcommand{\nn}{\mathbb{N}}
\newcommand{\zz}{\mathbb{Z}}

\newcommand{\qqand}{\qquad\text{and}\qquad}

\newcommand{\adjr}[5][30pt]{%
\begin{tikzcd}[ampersand replacement=\&, column sep=#1]
  #2\ar[r, shift left=7pt, "#3"]\&
  #5\ar[l, shift left=7pt, "#4"]
  \ar[l, phantom, "\scriptstyle\Leftarrow"]
\end{tikzcd}
}

\newenvironment{idris}
  {\VerbatimEnvironment
  	\begin{minted}[escapeinside=??, mathescape=true,frame=single, framesep=5pt, tabsize=1
	]{Idris}}
  {\end{minted}}
  
\newcommand{\idr}[1]{\mintinline{Idris}{#1}}

\setsecheadstyle{\bfseries\large\raggedright}
\setsubsecheadstyle{\bfseries\raggedright}
\firmlists
\date{\vspace{-3ex}}

\begin{document}
\title{$\poly$: An abundant categorical setting\\for mode-dependent dynamics}
\author{David I.\ Spivak}
\pagestyle{companion}

\maketitle

\begin{abstract}
Dynamical systems---by which we mean machines that take time-varying input, change their state, and produce output---can be wired together to form more complex systems. Previous work has shown how to allow collections of machines to reconfigure their wiring diagram dynamically, based on their collective state. This notion was called ``mode dependence'', and while the framework was compositional (forming an operad of re-wiring diagrams and algebra of mode-dependent dynamical systems on it), the formulation itself was more ``creative'' than it was natural.

In this paper we show that the theory of mode-dependent dynamical systems can be more naturally recast within the category $\poly$ of polynomial functors. This category is almost superlatively abundant in its structure: for example, it has \emph{four} interacting monoidal structures $(+,\times,\otimes,\circ)$, two of which ($\times,\otimes$) are monoidal closed, and the comonoids for $\circ$ are precisely categories in the usual sense. We discuss how the various structures in $\poly$ show up in the theory of dynamical systems. We also show that the usual coalgebraic formalism for dynamical systems takes place within $\poly$. Indeed one can see coalgebras as special dynamical systems---ones that do not record their history---formally analogous to contractible groupoids as special categories.
\end{abstract}

\chapter{Introduction}

We propose the category $\poly$ of polynomial functors on $\smset$ as a setting in which to model very general sorts of dynamics and interaction. Let's back up and say what exactly it is that we're generalizing.

A wiring diagram can be used to specify a fixed communication pattern between systems:
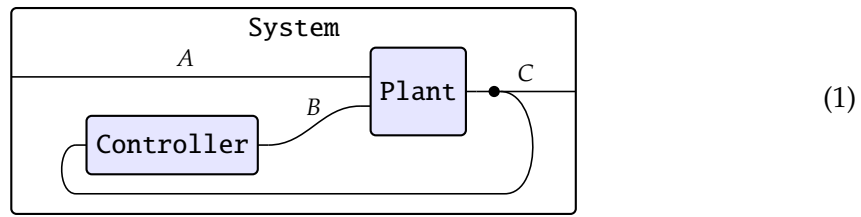
\begin{equation}\label{eqn.control_diag}
\begin{tikzpicture}[oriented WD, baseline=(B)]
	\node[bb={2}{1}, fill=blue!10] (plant) {\texttt{Plant}};
	\node[bb={1}{1}, below left=-1 and 1 of plant, fill=blue!10]  (cont) {\texttt{Controller}};
	\node[circle, inner sep=1.5pt, fill=black, right=.1] at (plant_out1) (pdot) {};
	\node[bb={0}{0}, inner ysep=15pt, inner xsep=1cm, fit=(plant) (pdot) (cont)] (outer) {};
	\coordinate (outer_out1) at (outer.east|-plant_out1);
	\coordinate (outer_in1) at (outer.west|-plant_in1);
	\begin{scope}[above, font=\footnotesize]
  	\draw (outer_in1) -- node {$A$} (plant_in1);
  	\draw (cont_out1) to node (B) {$B$} (plant_in2);
  	\draw (plant_out1) to node {$C$} (outer_out1);
  	\draw
  		let 
  			\p1 = (cont.south west-| pdot),
  			\p2 = (cont.south west),
  			\n1 = \bby,
  			\n2 = \bbportlen
  		in
  			(pdot) to[out=0, in=0]
  			(\x1+\n2, \y1-\n1) --
  			(\x2-\n2, \y2-\n1) to[out=180, in=180]
  			(cont_in1);
		\end{scope}
	\node[below=0of outer.north] {\texttt{System}};
\end{tikzpicture}
\end{equation}
Shown here, the plant---say a power plant or a car---is a dynamical system that receives input of type $A$ from the outside world and input of type $B$ from the controller, and it produces output of type $C$; this in turn is fed both to the outside world and to the controller. Given these fixed sets $A,B,C$, we will see shortly that the two interior boxes and one exterior box shown in \eqref{eqn.control_diag} can be faithfully represented by polynomials in one variable $\yon$, as follows:
\begin{equation}\label{eqn.basic_diagram}
	\mathtt{Plant}=C\yon^{AB}
	\qquad\quad
	\mathtt{Controller} = B\yon^C
	\qquad\quad
	\mathtt{System} = C\yon^A.
\end{equation}
Observe that in each case the output type is the coefficient on $\yon$, and the input type is the exponent on $\yon$. In Section~\ref{sec.wd_mode_dep} we will see that the wiring diagram \eqref{eqn.control_diag} itself, as well as the interacting dynamics, can be represented by morphisms involving these polynomials.

\section{Introduction to mode-dependence}

Notice that the polynomials in \eqref{eqn.basic_diagram} are monomials; it is this we want to generalize. By using more general polynomials such as $\mathtt{Robot}=\yon^{A_1A_2}+\yon+B\yon$, we can create a system for which the input-output types are not fixed: 
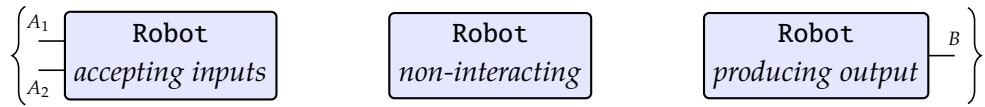
\[
\begin{tikzpicture}[oriented WD, align=center, bb port length=10pt, fill=blue!10]
	\node[bb={2}{0}, fill=blue!10] (in) {\texttt{Robot}\\\emph{accepting inputs}};
	\node[bb={0}{0}, right = of in, fill=blue!10] (busy) {\texttt{Robot}\\\emph{non-interacting}};
	\node[bb={0}{1}, right = of busy, fill=blue!10] (out) {\texttt{Robot}\\\emph{producing output}};
	\draw [decorate,decoration={brace,amplitude=5pt}]
	($(in.south west)+(-5pt-\bbportlen,-2pt)$) -- ($(in.north west)+(-5pt-\bbportlen,2pt)$);
	\draw [decorate,decoration={brace,amplitude=5pt}]
	($(out.north east)+(5pt+\bbportlen,2pt)$) -- ($(out.south east)+(5pt+\bbportlen,-2pt)$);
	\node[above=0 of in_in1, font=\tiny] {$A_1$};
	\node[below=0 of in_in2, font=\tiny] {$A_2$};
	\node[above=0 of out_out1, font=\tiny] {$B$};
\end{tikzpicture}
\]
What we discuss in this paper are dynamical systems whose interfaces change in time, and similarly where the wiring diagram connecting the systems changes in time. These changes will be based on the internal states of the systems---say robots---involved.

The real world is filled with instances of systems with time-varying input-output patterns. 
The network topology---the way that the system wires up---changes based on both internal and environmental contexts. Consider the following situations:\label{page.six}
\begin{enumerate}[itemsep=0pt]
	\item When too much force is applied to a material, bonds can break;
\end{enumerate}
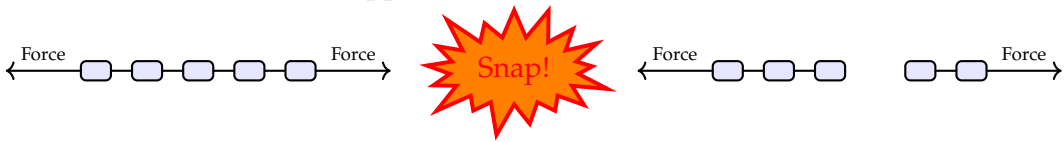
\[
\begin{tikzpicture}[oriented WD, bb small, bb port length=0]
	\foreach \i in {0,...,4} {
		\node[bb={1}{1}, fill=blue!10] at (1.7*\i,0) (X\i) {};
	}
	\foreach \i in {0,...,3} {
		\draw[thick] (X\i_out1) -- (X\the\numexpr\i+1\relax_in1);
	};
	\draw[thick, ->] (X0_in1) -- node[above, font=\tiny] {Force} +(-2.5,0);
	\draw[thick, ->] (X4_out1) -- node[above, font=\tiny] {Force} +(2.5,0) node (R) {};
\def\x{21};
	\foreach \i in {0,...,2} {
		\node[bb={1}{1}, fill=blue!10] at (\x+1.7*\i,0) (Y\i) {};
	}
	\foreach \i in {3,...,4} {
		\node[bb={1}{1}, fill=blue!10] at (\x+1.3+1.7*\i,0) (Y\i) {};
	}
	\foreach \i in {0,1,3} {
		\draw[thick] (Y\i_out1) -- (Y\the\numexpr\i+1\relax_in1);
	};
	\draw[thick, ->] (Y0_in1) -- node[above, font=\tiny] {Force} +(-2.5,0) node (L) {};
	\draw[thick, ->] (Y4_out1) -- node[above, font=\tiny] {Force} +(2.5,0);
	\node[starburst, draw, minimum width=2cm, minimum height=1.5cm,red,fill=orange,line width=1.5pt] at ($(L)!.5!(R)$)
{Snap!};
\end{tikzpicture}
\]
\begin{enumerate}[resume]
	\item A company may change its supplier at any time;
\end{enumerate}
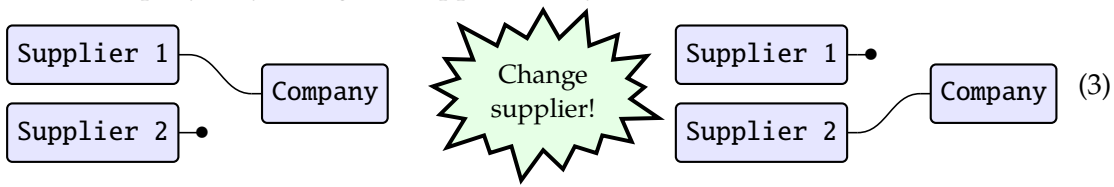
\begin{equation}\label{eqn.supplier}
\begin{tikzpicture}[oriented WD, font=\ttfamily, every node/.style={fill=blue!10}, baseline=(c)]
	\node[bb={0}{1}] (s1) {Supplier 1};
	\node[bb={0}{1}, below=of s1] (s2) {Supplier 2};
	\coordinate (helper) at ($(s1)!.5!(s2)$);
	\node[bb={1}{0}, right=1.5 of helper] (c) {Company};
	\draw (s1_out1) to (c_in1);
	\draw (s2_out1) to +(5pt,0) node[fill=none] {$\bullet$};
\begin{scope}[xshift=3.5in]
	\node[bb={0}{1}] (s1') {Supplier 1};
	\node[bb={0}{1}, below=of s1'] (s2') {Supplier 2};
	\coordinate (helper') at ($(s1')!.5!(s2')$);
	\node[bb={1}{0}, right=1.5 of helper'] (c') {Company};
	\draw (s2'_out1) to (c'_in1);
	\draw (s1'_out1) to +(5pt,0) node[fill=none] {$\bullet$};
\end{scope}
	\node[starburst, draw, minimum width=2cm, minimum height=2cm,align=center,fill=green!10, font=\small,line width=1.5pt] at ($(c)!.5!(helper')$)
{Change\\supplier!};
\end{tikzpicture}
\end{equation}
\begin{enumerate}[resume]
	\item When someone assembles a machine, their own outputs dictate the connection pattern of the machine's components.
\end{enumerate}
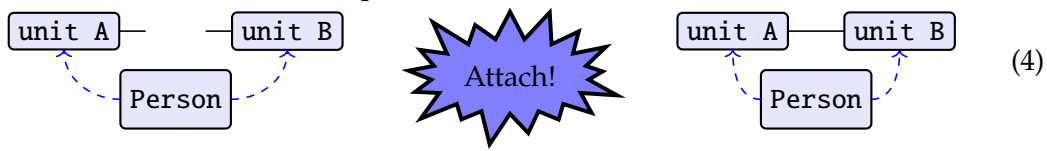
\begin{equation}\label{eqn.someone}
\begin{tikzpicture}[oriented WD, font=\ttfamily, bb port length=0, every node/.style={fill=blue!10}, baseline=(someone.north)]
	\node[bb port sep=.5, bb={0}{1}] (A) {unit A};
	\node[bb port sep=.5, bb={1}{0}, right=of A] (B) {unit B};
	\coordinate (helper) at ($(A)!.5!(B)$);
	\node[bb={1}{1}, below=2 of helper] (someone) {Person};
	\draw[->, dashed, blue] (someone_in1) to[out=180, in=270] (A.270);
	\draw[->, dashed, blue] (someone_out1) to[out=0, in=270] (B.270);
	\draw (A_out1) -- +(10pt,0);
	\draw (B_in1) -- +(-10pt,0);
\begin{scope}[xshift=3.5in]
	\node[bb port sep=.5, bb={0}{1}] (A') {unit A};
	\node[bb port sep=.5, bb={1}{0}, right=.5of A'] (B') {unit B};
	\coordinate (helper') at ($(A')!.5!(B')$);
	\node[bb={1}{1}, below=2 of helper'] (someone') {Person};
	\draw[->, dashed, blue] (someone'_in1) to[out=180, in=270] (A'.270);
	\draw[->, dashed, blue] (someone'_out1) to[out=0, in=270] (B'.270);
	\draw (A'_out1) -- (B'_in1);
\end{scope}
	\node[starburst, draw, minimum width=2cm, minimum height=2cm,fill=blue!50,line width=1.5pt, align=center, font=\upshape] at ($(B)!.5!(A')-(0,.6cm)$)
{Attach!};
\end{tikzpicture}
\end{equation}
We will discuss \eqref{eqn.supplier} and \eqref{eqn.someone} further in Example~\ref{ex.mode_dep}. In each of the above cases the wiring diagram---the connection pattern---changes based on the states (position, decision-making, environmental context, etc.) of some or all the systems involved. In \cite{spivak2017nesting} this was called \emph{mode-dependence}; the goal of that article was to create an operadic framework in which mode-dependent dynamics and communication could be specified compositionally. While successful, the presentation was fairly ad hoc. The purpose of the present paper is to explain that the category $\poly$ provides an abundant setting in which to work quite naturally with mode-dependent dynamics.

When we say that $\poly$ is abundant, we mean that it is exceptionally rich in structure, and that structure is highly relevant to dynamical systems. Here are some of the features of this category:
\begin{enumerate}[itemsep=0pt]
	\item $\poly$ has coproducts and products, $+, \times$, the usual sum and product of polynomials.
	\item $\poly$ has two additional monoidal structures: $\otimes$ and $\circ$. 
	\item $\poly$ has two monoidal closed structures: for $\times$ (cartesian closure) and $\otimes$.
	\item $\poly$ has a duoidal structure: $(\circ)\otimes(\circ)\to(\otimes)\circ(\otimes)$.
	\item $\poly$ has all small limits and is extensive.
	\item $\poly$ has two orthogonal factorization systems (epi/mono and vertical/cartesian).
	\item $\poly$ admits a monoidal bifibration $\poly\to\smset$ with $\otimes\mapsto\times$.
	\item $\poly$ admits an adjoint quadruple with $\smset$ and an adjoint pair with $\smset\op$.
	\item Comonoids in $(\poly,\circ)$ are precisely categories in the usual sense.
\end{enumerate}

In Section~\ref{chap.poly_intro} we will introduce $\poly$ and many of its interesting features. In Section~\ref{chap.poly_mode_dep_dyn}, we will discuss how these features relate to dynamical systems.

\section{Acknowledgments}
We thank Richard Garner and David Jaz Myers for helpful and interesting conversations. We also appreciate support from Honeywell and AFOSR grants FA9550-17-1-0058 and FA9550-19-1-0113.

\chapter{Introduction to $\poly$}\label{chap.poly_intro}

\section{Polynomial functors}

\nocite{girard1988normal}\nocite{moerdijk2000wellfounded}

\begin{notation}
We usually denote sets with upper-case letters $A,B$, etc.; the exception is ordinals: we denote the $n$th ordinal by $n=\{\fun{1},\ldots,\fun{n}\}$. We denote functions between sets---including elements of sets---using upright letters $\fun{f}\colon A\to B$ and $\fun{a}\in A$. 

All polynomials discussed here have a single variable, always $\yon$; in particular $\yon$ itself is a polynomial. Coefficients and exponents of polynomials are arbitrary sets, e.g.\ $\nn\yon^\rr+3$ is a polynomial. Every set $A$ will also be a polynomial, namely a constant. We denote generic polynomials with lower-case letters $p,q$, etc. 
\end{notation}

Recall that a \emph{representable functor} $\smset\to\smset$ is one of the form $\smset(A, -)$ for a set $A$.
We denote this functor by $\yon^A\colon\smset\to\smset$ and say it is \emph{represented by} $A\in\smset$. For example $\yon^3$ is represented by $3$ and $\yon^3(2)\cong8$. As $A$ varies we obtain the contravariant Yoneda embedding.

Classically, a polynomial $p$ in one variable with set coefficients is a function $p(\yon)=A_n\yon^n+\cdots+A_1\yon^1+A_0\yon^0$ with each $A_i\in\nn$. In category theory this is often generalized to allow for infinitely many terms and infinite exponents; e.g.\ we consider the following to be a polynomial
\begin{equation}\label{eqn.poly}
p(\yon)=\sum_{\fun{i}\in I}\yon^{A_\fun{i}}
\end{equation}
for arbitrary small sets $I$ and $A$. We can think of such a $p$ as a functor $\smset\to\smset$; it sends a set $X\in\ob(\smset)$ to the coproduct, over $\fun{i}\in I$, of the set $X^{A_\fun{i}}$ of functions $A_\fun{i}\to X$, or equivalently the $A_\fun{i}$-fold product of $X$ with itself. The result is covariantly functorial in $X$.

Considered this way, $p$ is called a \emph{polynomial functor}; polynomial functors sit inside of the category of all functors $\smset\to\smset$ as a full subcategory, namely the one spanned by coproducts of representables.

\begin{definition}
The category $\poly$ has polynomial functors $p(\yon)$ as in \eqref{eqn.poly} as objects and natural transformations between them as morphisms.
\end{definition}
In $\poly$, products distribute over coproducts, $\Big(\sum_i p_i\Big)\times q  \cong \sum_i\left(p_i\times q\right)$. More generally, for any discrete category $A$, functor $I\colon A\to\smset$, and functor $p\colon\sum_{a\in A}I(a)\to\poly$ (the above specific case being $A\coloneqq 2$, $I(2)\coloneqq 1$, $p_{1,i}\coloneqq p_i$, and $p_{2,1}\coloneqq q$) there is an isomorphism
\begin{equation}\label{eqn.complete_dist}
  \prod_{a\in A}\;\;\sum_{i\in I(a)}\;p_{(a,i)}
  \quad\cong\quad
  \sum_{i\in\prod_{a\in A}I(a)}\;\;\prod_{a\in A}\;p_{(a,i(a))}.
\end{equation}
In fact, $\poly$ can be characterized as the free category that has coproducts, products, and satisfies Eq.~\ref{eqn.complete_dist}, (called \emph{complete distributivity} or the Axiom of Choice; see e.g.\ \cite{jacobs1999categorical}). $\poly$ is also equivalent to the Grothendieck construction of the canonical functor $\smset\op\to\smcat$ sending each object to the corresponding slice category (opposite) $A\mapsto(\smset/A)\op$ and sending $\fun{f}\colon B\to A$ to pullback along $\fun{f}$.

\begin{notation}
We denote the product of polynomials by juxtaposition or sometimes $\mdot$, i.e.\
\[pq\coloneqq p\times q=p\mdot q.\]
For any set $A$ we denote the $A$-fold repeated product of $p$ by $p^A\coloneqq \prod_{\fun{a}\in A}p$; in particular $p^1\cong p$ and $p^0\cong 1$. The representable $\yon^n$ is indeed the $n$-fold repeated product of $\yon$.
\end{notation}

For any polynomial $p$, the set $p(1)$  has particular importance; it can be identified with the set of representable summands (pure-power terms) $\yon^k$ in $p$. For example if $p=\yon^2+3\yon+2$ then $p(1)=6$ corresponding to the six representable summands in $p=\yon^2+\yon^1+\yon^1+\yon^1+\yon^0+\yon^0$. We will denote the representing object for the $\fun{i}$th representable summand of $p$ by $\fun{p_i}$, i.e.\
\[
p=\sum_{\fun{i}\in p(1)}\yon^{\fun{p_i}}.
\]

There are many ways to think about polynomials, and one becomes more versatile by being able to use different representations for different purposes. So far we have been writing polynomials in the typical algebraic style, but one can also represent them as bundles, as forests of corollas, or as dependent types. 
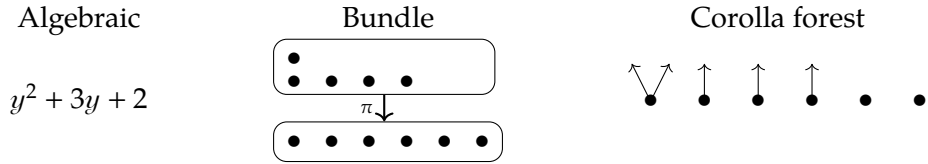
\[
\begin{array}{ccccc}
\text{Algebraic}&~\qquad~&\text{Bundle}&~\qquad~&\text{Corolla forest}\\
\yon^2+3\yon+2&&
\begin{tikzpicture}[baseline=(pi)]
  \foreach \n/\i in {1/2,2/1,3/1,4/1,5/0,6/0} {
  	\node at (.5*\n,0) (b\n) {$\bullet$};
		\node (b\n0) at (.5*\n,.8) {};
  	\ifthenelse{\i>0}{
		\foreach \j in {1,...,\i} {
		\node (b\n\j) at (.5*\n,.5+.3*\j) {$\bullet$};
		};
		}{};
  };
  \node[draw, inner sep=1pt, rounded corners=5, fit=(b1) (b6)] (B)  {};
  \node[draw, inner sep=1pt, rounded corners=5, fit=(b12) (b60)] (E) {};
  \draw[->, thick] (E) -- node[left, font=\tiny] (pi) {$\pi$} (B.north-|E);
\end{tikzpicture}
&&
\begin{tikzpicture}[trees, grow'=up]
  \node (1) {$\bullet$} 
    child {}
    child {};
  \node[right=.5 of 1] (2) {$\bullet$} 
    child {};
  \node[right=.5 of 2] (3) {$\bullet$} 
    child {};
  \node[right=.5 of 3] (4) {$\bullet$} 
    child {};
  \node[right=.5 of 4] (5) {$\bullet$};
  \node[right=.5 of 5] (6) {$\bullet$};
\end{tikzpicture}
\end{array}
\]
Given a bundle $\pi\colon E\to B$, and element $\fun{b}\in B$, we denote the fiber $\pi\inv(\fun{b})$ by $E_\fun{b}$. We will refer to elements of $B$ as \emph{positions} and elements of $E_\fun{b}$ as the \emph{directions} in position $\fun{b}$. From the algebraic viewpoint, a position is a `pure-power', or representable summand, and the associated direction-type is its `exponent' or representing object; from the tree viewpoint, a position is a \emph{root} and the associated directions are its \emph{leaves}.

Polynomials can be implemented in a dependently typed programming language, such as Idris. Here is a specification of the type for polynomials:%
\begin{idris}
record Poly where                     
   constructor MkPoly              -- To construct a poly, define:
   position  : Type                -- the "positions" (as a type), and
   direction : position -> Type    -- the "directions" in each position.
\end{idris}
For example, the expression \idr{MkPoly Integer (\i => Double)} means that the type of positions is $\zz$ and for each position the type of directions is the type \idr{double} of double-precision floating point numbers; thinking of it as the reals, this denotes the monomial $\zz\yon^\rr$. We will only discuss Idris once more in this document, though almost everything we discuss has been implemented; please write to the author for more information.

\section{Morphisms of polynomials, concretely}
As mentioned, the morphisms between polynomials are the natural transformations. As easy as this is to state---and as much as it gives us confidence in the reasonableness of the definition---it can be useful to have a more hands-on understanding of the morphisms.

By the Yoneda lemma, a morphism $\yon^A\to\yon^B$ can be identified with a function $B\to A$. One can prove that $+$ is the coproduct in $\poly$ and $\times$ is the product. Thus $\yon^2+3\yon+2$ is a product of $\yon+1$ and $\yon+2$, and it is a coproduct of $\yon^2+1$ and $3\yon+1$. This also holds for infinite sums and products: the usual algebraic operations coincide with the categorical operations. From this, and the fact that coproducts of functors $\smset\to\smset$ are taken pointwise, we obtain the following formula for the set of morphisms $p\to q$:
\[\poly(p,q)\cong\prod_{\fun{i}\in p(1)}\;\sum_{\fun{j}\in q(1)}\fun{p_i}^{\fun{q_j}}.\]
For each representable summand of $p$---i.e.\ position of $p$---choose a representable summand of $q$ and give a function from the representing object (exponent) in $q$ back to the representing object (exponent) in $p$. Thus for example $\poly(\yon^2+3\yon+2,\ \yon^5+1)\cong(2^5+1)(1^5+1)(1^5+1)(1^5+1)(0^5+1)(0^5+1)$.

In terms of bundles, a morphism $E\to B$ to $E'\to B'$ consists of a pair $(\fun{f},\fun{f}^\sharp)$ as shown:
\[
\begin{tikzcd}[row sep=15pt]
  E\ar[d]&
  B\times_{B'}E'\ar[r]\ar[d]\ar[l, "\fun{f}^\sharp"']&
  E'\ar[d]\\
  B\ar[r, equal]&
  B\ar[r, "\fun{f}"']&
  B'\ar[ul, phantom, very near end, "\lrcorner"]
\end{tikzcd}
\]
This will be the most convenient way to write morphisms of polynomials; we further denote by $\fun{f}^\sharp_\fun{i}$ the map on fibers $E'(f(\fun{p}))\to E(\fun{p})$. We refer to $\fun{f}$ as the \emph{on-positions} function and $\fun{f}^\sharp$ as the \emph{on-directions} function. This way of thinking about morphisms of polynomials extends readily to Idris:
\begin{idris}
record Lens (dom : Poly) (cod : Poly) where
   constructor MkLens
   onPos : position dom -> position cod
   onDir : (i : position dom) -> direction cod j -> direction dom i
      where j = onPos i
\end{idris}

The reason for the name \emph{lens} comes from the following; see also \cites{abbott2003categories,spivak2019generalized}.
\begin{example}[Bimorphic lenses]
In \cite{hedges2018limits}, Hedges defines the category of bimorphic lenses to have objects given by pairs of sets $(A,B)$ and morphisms (called \emph{lenses}) from $(A,B)$ to $(A',B')$ defined by a pair of maps $A\to A'$ and $A\times B'\to B$. It is straightforward to check that Hedges' category of bimorphic lenses is equivalent to the full subcategory of $\poly$ spanned by the monomials $B\yon^A$.

  Monomials $B\yon^A$ will play a special role in the theory of this paper, namely they correspond to interfaces that have \emph{fixed inputs ($A$) and outputs ($B$)}, e.g.\ as seen in \eqref{eqn.control_diag}.
\end{example}

The category $\poly$ has all small limits. Suppose given a small category $J$ and functor $p\colon J\to\poly$, and for each $\fun{j}$, let $p^\fun{j}$ denote the corresponding polynomial. The limit $\lim_{\fun{j}\in J} p$ has positions given by the limit $\lim_{\fun{j}\in J}p^\fun{j}(1)$ of positions, and for each such position $(\fun{i^j})_{\fun{j}\in J}$, where $\fun{i^j}\in p^\fun{j}(1)$, the set of directions there is given by the colimit $\colim_{\fun{j}\in J}\fun{p^j_{i^j}}$ of directions.

We note two orthogonal factorization systems on $\poly$: (epi/mono) and (vertical/cartesian). The first is straightforward (e.g.\ epimorphisms of polynomials are surjective on positions and injective on directions). More interestingly, the functor $p\mapsto p(1)$ is a monoidal *-bifibration in the sense of \cite[Definition 12.1]{shulman2008framed}. Indeed, if $B=p(1)$ and we have a function $f\colon A\to B$, we can take the pullback of polynomials
\[
\begin{tikzcd}[row sep=15pt]
	A\times_Bp\ar[r, "\mathrm{cart}_f"]\ar[d]&
	p\ar[d]\\
	A\ar[r, "f"']&
	B\ar[ul, phantom, very near end, "\lrcorner"]
\end{tikzcd}
\]
Thus we obtain a fibration $\poly\to\smset$, with its attendant vertical/cartesian factorization system. Moreover each functor $f^*\colon\poly_B\to\poly_A$ has both a left adjoint $f_!$ and a right adjoint $f_*$, and both $f_!$ and $f^*$ interact well with $\otimes$, a monoidal product we will introduce in Section~\ref{sec.monoidal_structures}. In fact, identifying $\poly_A\op$ with $\smset^A$, the functors $\smset^I\to\smset^J$ arising from multivariate polynomials $I\From{f} E\To{g} B\To{h}J$ as in \cite{kock2012polynomial} can be represented using the $*$-bifibration structure, namely as $(h_*g_!f^*)\op$.

\section{Adjunctions with $\smset$ and $\smset\op$}

It is useful to note that $\poly$ contains two copies of $\smset$ and a copy of $\smset\op$, namely as the constant polynomials $A$, the linear polynomials $A\yon$, and the representables $\yon^A$. Indeed there is an adjoint quadruple and an adjoint pair as follows, labeled by where they send objects $A\in\smset$, $p\in\poly$:%
\footnote{
	We use the notation 
	$\begin{tikzcd}[ampersand replacement=\&]
		C\ar[r, shift left=4pt, "L"]\&
		D\ar[l, shift left=4pt, "R"]\ar[l, phantom, "\scriptstyle\Rightarrow"]
	\end{tikzcd}$
	to denote an adjunction $L\dashv R$. The double arrow, always pointing in the direction of the left adjoint, indicates both the unit $C\Rightarrow R\circ L$ and the counit $L\circ R\Rightarrow D$ of the adjunction.
}
\begin{equation}\label{eqn.adjoints_galore}
\begin{tikzcd}[column sep=60pt]
  \smset
  	\ar[r, shift left=7pt, "A" description]
		\ar[r, shift left=-21pt, "A\yon"']&
  \poly
  	\ar[l, shift right=21pt, "p(0)"']
  	\ar[l, shift right=-7pt, "p(1)" description]
	\ar[l, phantom, "\scriptstyle\Leftarrow"]
	\ar[l, phantom, shift left=14pt, "\scriptstyle\Rightarrow"]
	\ar[l, phantom, shift right=14pt, "\scriptstyle\Rightarrow"]
\end{tikzcd}
\hspace{.6in}
\adjr[50pt]{\smset\op}{\yon^A}{\Gamma p}{\poly}.
\end{equation}
All of the functors out of $\smset$ and $\smset\op$ shown in \eqref{eqn.adjoints_galore} are fully faithful, and the rightmost adjoint $p\mapsto p(0)$ preserves coproducts. The functor $\Gamma$ is given by \emph{global sections}: $\Gamma p\coloneqq\poly(p,\yon)\cong\prod_{\fun{i}\in p(1)}\fun{p_i}$.

For each $A\in\smset$ the functor $\poly\to\smset$ given by $q\mapsto q(A)$ has a left adjoint, namely $B\mapsto B\yon^A$; we saw this for the cases $A\cong 0,1$ in Eq.~\ref{eqn.adjoints_galore}. Using $p\coloneqq\yon^A$ and the Yoneda lemma, this generalizes to a two-variable adjunction $\smset\times\poly\to\poly$:
\begin{align}\label{eqn.two_var_adj}
\poly(Ap,q)\cong\poly(p,q^A)\cong\smset(A,\poly(p,q)).
\end{align}

\section{Monoidal structures on $\poly$}\label{sec.monoidal_structures}

We have already mentioned two monoidal structures on $\poly$, namely coproduct $(+,0)$ and product $(\times,1)$. They are given by the following formulas:
\begin{equation}\label{eqn.sums_prods}
  p+ q=\sum_{\fun{i}\in p(1)}\yon^{\fun{p_i}}+\sum_{\fun{j}\in q(1)}\yon^{\fun{q_j}}
  \qqand
  p\times q=\sum_{\fun{i}\in p(1)}\sum_{\fun{j}\in q(1)}\yon^{\fun{p_i+q_j}}.
\end{equation}
These form a distributive category. The product monoidal structure is closed---$\poly$ is cartesian closed---and we denote this closure operation by exponentiation:
\begin{equation}\label{eqn.exponentiate}
q^p\cong\prod_{\fun{i}\in p(1)}q\circ (\fun{p_i}+\yon).
\end{equation}
Thus for example $(\yon^2+3\yon+2)^{\yon^5+\yon^4}\cong\big((5+\yon)^2+3(5+\yon)+2\big)\mdot\big((4+\yon)^2+3(4+\yon)+2\big)$. The constant-polynomials functor $\smset\to\poly$ is cartesian closed.

In terms of bundles, the coproduct is given by disjoint union, and product is given by adding fibers (though the formula is reminiscent of adding fractions):
\[
  \left(
  \begin{tikzcd}
  E\ar[d]\\B
  \end{tikzcd}
  \right)
+
  \left(
  \begin{tikzcd}
  E'\ar[d]\\B'
  \end{tikzcd}
  \right)
\cong
  \left(
  \begin{tikzcd}
  E+ E'\ar[d]\\B+ B'
  \end{tikzcd}
  \right)
\hspace{1in}
  \left(
  \begin{tikzcd}
  E\ar[d]\\B
  \end{tikzcd}
  \right)
\times
  \left(
  \begin{tikzcd}
  E'\ar[d]\\B'
  \end{tikzcd}
  \right)
\cong
  \left(
  \begin{tikzcd}
  E\times B'+B\times E'\ar[d]\\B\times B'
  \end{tikzcd}
  \right)
\]
In terms of forests, coproduct (undrawn) is given by disjoint union and product is given by multiplying the roots and adding the leaves. Here is a picture of $(\yon+1)(\yon+2)\cong\yon^2+3\yon+2$:
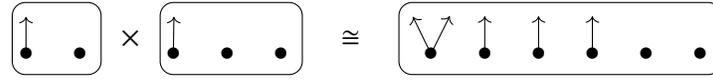
\[
\begin{tikzpicture}
	\node (p1) {
	\begin{tikzpicture}[trees, grow'=up]
  	\node (1) {$\bullet$}
  		child {node (a) {}};
  	\node[right=.5 of 1] (2) {$\bullet$};
    \node[draw, rounded corners=5pt, inner sep=5pt, fit=(a) (2)] (x) {};	
	\end{tikzpicture}
	};
	\node (p2) [right=.5 of p1] {
	\begin{tikzpicture}[trees, grow'=up]
  	\node (1) {$\bullet$}
  		child {node (a) {}};
  	\node[right=.5 of 1] (2) {$\bullet$};
  	\node[right=.5 of 2] (3) {$\bullet$};
    \node[draw, rounded corners=5pt, inner sep=5pt, fit=(a) (3)] (x) {};	
	\end{tikzpicture}
	};
	\node (p3) [right=1 of p2] {
	\begin{tikzpicture}[trees, grow'=up]
  	\node (1) {$\bullet$} 
      child {node (1a) {}}
      child {};
    \node[right=.5 of 1] (2) {$\bullet$} 
      child {};
    \node[right=.5 of 2] (3) {$\bullet$} 
      child {};
    \node[right=.5 of 3] (4) {$\bullet$} 
      child {};
    \node[right=.5 of 4] (5) {$\bullet$};
    \node[right=.5 of 5] (6) {$\bullet$};
    \node[draw, rounded corners=5pt, inner sep=5pt, fit=(1a) (6)] (prod) {};
	\end{tikzpicture}
   };
   \node at ($(p1.east)!.5!(p2.west)$) {$\times$};
   \node at ($(p2.east)!.5!(p3.west)$) {$\cong$};
\end{tikzpicture}
\]

There are two more monoidal structures on $\poly$; one is symmetric and is denoted $(\otimes,\yon)$, and the other is not symmetric and is denoted $(\circ,\yon)$. We first discuss $\otimes$. In terms of polynomials, it is given by the \emph{Dirichlet product}%
\footnote{The reason for the name \emph{Dirichlet} is that if one replaces polynomials with Dirichlet series by reversing each summand $\yon^A$ to $A^\yon$, the result is the usual product. For example
\[
(3^\yon+2^\yon)\times(4^\yon+0^\yon)\cong12^\yon+8^\yon+2\mdot0^\yon
\]
See \cite{spivak2020dirichlet} for more on the connection between Dirichlet series and polynomials.
}
\begin{equation}\label{eqn.dirichlet}
p\otimes q=\sum_{\fun{i}\in p(1)}\sum_{\fun{j}\in q(1)}\yon^{\fun{p_iq_j}},
\end{equation}
which we invite the reader to compare with $\times$ from Eq.~\ref{eqn.sums_prods}. For example $(\yon^3+\yon)\otimes(y^2+\yon^0)\cong\yon^{6}+\yon^2+2\yon^0$. Like $\times$, the Dirichlet product $\otimes$ distributes over $+$. In terms of bundles, Dirichlet product is straightforward:
\[
  \left(
  \begin{tikzcd}
  E\ar[d]\\B
  \end{tikzcd}
  \right)
\otimes
  \left(
  \begin{tikzcd}
  E'\ar[d]\\B'
  \end{tikzcd}
  \right)
\cong
  \left(
  \begin{tikzcd}
  E\times E'\ar[d]\\B\times B'
  \end{tikzcd}
  \right).
\]
In terms of forests, one multiplies roots and for each pair, multiplies the leaves:
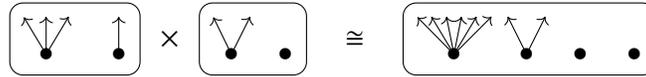
\[
\begin{tikzpicture}
	\node (p1) {
	\begin{tikzpicture}[trees, grow'=up]
  	\node (1) {$\bullet$} [sibling distance=.3cm] 
  		child {node (a) {}}
			child
			child;
  	\node[right=.75 of 1] (2) {$\bullet$}
			child;
    \node[draw, rounded corners=5pt, inner sep=5pt, fit=(a) (2)] (x) {};	
	\end{tikzpicture}
	};
	\node (p2) [right=.5 of p1] {
	\begin{tikzpicture}[trees, grow'=up]
  	\node (1) {$\bullet$}
  		child {node (a) {}}
			child;
  	\node[right=.5 of 1] (2) {$\bullet$};
    \node[draw, rounded corners=5pt, inner sep=5pt, fit=(a) (2)] (x) {};	
	\end{tikzpicture}
	};
	\node (p3) [right=1 of p2] {
	\begin{tikzpicture}[trees, grow'=up]
  	\node (1) {$\bullet$} [sibling distance=.2cm] 
      child {node (1a) {}}
      child
      child
      child
      child
      child {};
    \node[right=.75 of 1] (2) {$\bullet$} 
      child {}
      child;
    \node[right=.5 of 2] (3) {$\bullet$}; 
    \node[right=.5 of 3] (4) {$\bullet$};
    \node[draw, rounded corners=5pt, inner sep=5pt, fit=(1a) (4)] (prod) {};
	\end{tikzpicture}
   };
   \node at ($(p1.east)!.5!(p2.west)$) {$\times$};
   \node at ($(p2.east)!.5!(p3.west)$) {$\cong$};
\end{tikzpicture}
\]
The Dirichlet monoidal structure is closed as well and its formula is similar to that in \eqref{eqn.exponentiate}. We denote this closure operation (internal hom) using brackets:
\begin{equation}\label{eqn.bracket}
[p,q]\cong\prod_{\fun{i}\in p(1)}q\circ (\fun{p_i}\yon).
\end{equation}
Thus for example $[\yon^5+\yon^4,\yon^2+3\yon+2)]\cong((5\yon)^2+3(5\yon)+2)\mdot((4\yon)^2+3(4\yon)+2)$.

The last monoidal structure we discuss, $(\circ,\yon)$, was already used above in Eqs.~\ref{eqn.exponentiate}~and~\ref{eqn.bracket}. It is the usual composition of polynomials, both algebraically and as functors; e.g.\ $(\yon^2+\yon)\circ(\yon^3+1)\cong\yon^6+3\yon^3+2$. Thinking of $p$ as a functor, its evaluation at a set $A$ is $p\circ A$.

The most computationally useful formula for $p\circ q$ is probably the following:
\begin{equation}\label{eqn.strategy}
p\circ q\cong\sum_{\fun{i}\in p(1)}\;\prod_{\fun{d}\in\fun{p_i}}\;\sum_{\fun{j}\in q(1)}\;\prod_{\fun{e}\in\fun{q_j}}\;\yon.
\end{equation}
In terms of forests, $p\circ q$ is obtained by adding up all ways to adjoin trees in $q$ to leaves in $p$. For example, here is $(\yon^2+y)\circ(\yon^3+1)$:
\begin{equation}\label{eqn.trees_comp}
\begin{tikzpicture}[trees, grow'=up, baseline = (p1)]
	\node (p1) {
	\begin{tikzpicture}
  	\node (1) {$\bullet$}
  		child {node (a) {}}
			child {};
  	\node[right=.5 of 1] (2) {$\bullet$}
			child {};
    \node[draw, rounded corners=5pt, inner sep=5pt, fit=(a) (2)] (x) {};	
	\end{tikzpicture}
	};
	\node (p2) [right=.5 of p1] {
	\begin{tikzpicture}
  	\node (1) {$\bullet$}[sibling distance=.3cm]
  		child {node (a) {}}
			child {}
			child {};
  	\node[right=.5 of 1] (2) {$\bullet$};
    \node[draw, rounded corners=5pt, inner sep=5pt, fit=(a) (2)] (x) {};	
	\end{tikzpicture}
	};
	\node (p3) [right=1 of p2] {
  \begin{tikzpicture}[trees, grow'=up]
    \node (a) {$\bullet$}[sibling distance=.75cm] 
      child {[fill]
      	node[left=-3pt] {$\bullet$}[sibling distance=.3cm]
				  child {node (x) {}} child child 
			}
      child {[fill]
      	node[left=-3pt] {$\bullet$}[sibling distance=.3cm]
				  child child child
  		};
    \node (b) [right=1.5 of a] {$\bullet$} 
      child {[fill]
      	node[left=-3pt] {$\bullet$}[sibling distance=.3cm]
				  child child child
			}
      child {[fill]
      	node[left=-3pt] {$\bullet$}
  		};
    \node (c) [right=1 of b]{$\bullet$}
      child {[fill]
      	node[left=-3pt] {$\bullet$}
			}
      child {[fill]
      	node[left=-3pt] {$\bullet$}[sibling distance=.3cm]
				  child child child
  		};
    \node (d) [right=1 of c]{$\bullet$} 
      child {[fill]
      	node[left=-3pt] {$\bullet$}
			}
      child {[fill]
      	node[left=-3pt] {$\bullet$}
  		};		
    \node (e) [right=1 of d] {$\bullet$}
      child {[fill]
      	node[left=-3pt] {$\bullet$}[sibling distance=.3cm]
				  child child child 
			};
    \node (f) [right=1 of e] {$\bullet$}
      child {[fill]
      	node[left=-3pt] {$\bullet$}
			};
    \node[draw, fill=yellow!80!black, opacity=0.2, rounded corners=5pt, inner sep=5pt, fit=(x) (f)] (y) {};	
  \end{tikzpicture}
   };
   \node at ($(p1.east)!.5!(p2.west)$) {$\circ$};
   \node at ($(p2.east)!.5!(p3.west)$) {$\cong$};
\end{tikzpicture}
\end{equation}
More precisely, the monoidal operation $\circ$ collapses the trees in \eqref{eqn.trees_comp} to mere corollas:
\[
\begin{tikzpicture}[trees, grow'=up]
	\node (p1) {
	\begin{tikzpicture}
  	\node (1) {$\bullet$}
  		child {node (a) {}}
			child {};
  	\node[right=.5 of 1] (2) {$\bullet$}
			child {};
    \node[draw, rounded corners=5pt, inner sep=5pt, fit=(a) (2)] (x) {};	
	\end{tikzpicture}
	};
	\node (p2) [right=.5 of p1] {
	\begin{tikzpicture}
  	\node (1) {$\bullet$}[sibling distance=.3cm]
  		child {node (a) {}}
			child {}
			child {};
  	\node[right=.5 of 1] (2) {$\bullet$};
    \node[draw, rounded corners=5pt, inner sep=5pt, fit=(a) (2)] (x) {};	
	\end{tikzpicture}
	};
	\node (p3) [right=1 of p2] {
  \begin{tikzpicture}[trees, grow'=up]
    \node (a) {$\bullet$}[sibling distance=.2cm] 
      child {node (x) {}} child child child child child;
    \node (b) [right=1 of a] {$\bullet$} [sibling distance=.3cm]
      child child child;
    \node (c) [right=1 of b]{$\bullet$}[sibling distance=.3cm]
      child child child;
    \node (d) [right=1 of c]{$\bullet$};		
    \node (e) [right=1 of d] {$\bullet$}[sibling distance=.3cm]
      child child child;
    \node (f) [right=1 of e] {$\bullet$};
    \node[draw, rounded corners=5pt, inner sep=5pt, fit=(x) (f)] (y) {};	
  \end{tikzpicture}
   };
   \node at ($(p1.east)!.5!(p2.west)$) {$\circ$};
   \node at ($(p2.east)!.5!(p3.west)$) {$\cong$};
\end{tikzpicture}
\]

The composition product $\circ$ is \emph{duoidal} over $\otimes$ in the sense that there is a natural map
\begin{equation}\label{eqn.duoidal}
(p_1\circ p_2)\otimes(q_1\circ q_2)\to(p_1\otimes q_1)\circ(p_2\otimes q_2),
\end{equation}
satisfying the usual axioms. Both $+$ and $\times$ commute with $\circ$ on the left
\[(pq+r)\circ s\cong (p\circ s)(q\circ s)+(r\circ s).\]

\section{Comonoids for $\circ$ are categories}

A comonoid in the (nonsymmetric) monoidal category $(\poly,\circ,\yon)$ is a tuple $(C,\epsilon,\delta)$, where $C$ is a polynomial,%
\footnote{We use upper case to denote the polynomials that underlie comonoids.}
 and $\epsilon\colon C\to \yon$ and $\delta\colon C\to C\circ C$ are morphisms of polynomials, such that the usual diagrams commute. Using Eq.~\ref{eqn.duoidal}, we can lift the Dirichlet product on polynomials to a monoidal structure $(\otimes,\yon)$ on comonoids.

One of the most surprising aspects of $\poly$ is that the comonoids for $\circ$---polynomial comonads on $\smset$---are categories in the usual sense! This requires a calculation (see Theorem~\ref{thm.cats_comons}), though it can be visualized using tree composition \eqref{eqn.trees_comp}. Sums and Dirichlet products of comonoids correspond to coproducts and products of categories, respectively.

Note that morphisms of comonoids correspond not to functors but to \emph{cofunctors}, first defined in \cite{higgins1993duality}; see also \cite{aguiar1997internal}. This notion is not well-enough known, so we recall it. We temporarily use $C_0$ for the objects and $C_1$ for the morphisms in a category $C$.

\begin{definition}[Cofunctor]
Let $C$ and $D$ be categories. A \emph{cofunctor} $F\colon C\nrightarrow D$ consists of
\begin{enumerate}[itemsep=0pt]
  \item a function  $F_0\colon C_0\to D_0$ \emph{on objects} and
  \item a function $F^\sharp\colon C_0\times_{D_0}D_1\to C_1$ \emph{backwards on morphisms},
\end{enumerate}
satisfying the following conditions:%
\footnote{
The cofunctor laws written in diagram form are as follows:
\[
\begin{tikzcd}[column sep=15pt, ampersand replacement=\&]
  C_0\times_{D_0}D_0\ar[r, "\cong"]\ar[d, "\id_D"']\&
  C_0\ar[d, "\id_C"]\\
  C_0\times_{D_0}D_1\ar[r, "F^\sharp"']\&
  C_1\ar[ul, phantom, shift right=2pt, "\text{(i)}"]
\end{tikzcd}
\quad
\begin{tikzcd}[column sep=15pt, ampersand replacement=\&]
	C_0\times_{D_0}D_1\ar[r, "F^\sharp"]\ar[d, "\pi_2"']\&
	C_1\ar[r, "\cod"]\&
	C_0\ar[d, "F_0"]\\
	D_1\ar[rr, "\cod"']\ar[urr, phantom, "\text{(ii)}"]\&\&
	D_0
\end{tikzcd}
\quad
\begin{tikzcd}[column sep=15pt, ampersand replacement=\&]
	C_0\times_{D_0}D_1\times_{D_0}D_1\ar[r, "\circ_D"]\ar[d, "F^\sharp"']\&[-15pt]
	C_0\times_{D_0}D_1\ar[r, "F^\sharp"]\&
	C_1\\
	C_1\times_{D_0}D_1\ar[r, "\cong"']\&
	C_1\times_{C_0}C_0\times_{D_0}D_1\ar[r, "F^\sharp"']\ar[u, phantom, "\text{(iii)}"]\&
	C_1\times_{C_0}C_1\ar[u, "\circ_C"']
\end{tikzcd}
\]
}
\begin{enumerate}[itemsep=0pt, label=\roman*.]
	\item $F^\sharp(c,\id_{F_0c})=\id_c$ for any $c\in C_0$;
	\item $F_0\cod F^\sharp(c,g))=\cod g$ for any $c\in C_0$ and $g\colon F_0(c)\to\cod(g)$ in $D$;
	\item $F^\sharp(\cod F^\sharp(c,g_1),g_2)\circ F^\sharp(c,g_1)=F^\sharp(c,g_2\circ g_1)$ for composable arrows $g_1,g_2$ out of $F_0 c$.
\end{enumerate}
\end{definition}

\begin{theorem}[Ahman-Uustalu \cite{ahman2016directed}]\label{thm.cats_comons}
The following categories are equivalent:
\begin{enumerate}[itemsep=0pt]
	\item the category $\smcat^\sharp$ of categories and cofunctors;
	\item the category $\comon(\poly)$ of comonoids in $(\poly,\circ,\yon)$ and comonoid morphisms.
\end{enumerate}
\end{theorem}

The rough idea is that if $\sum_{i\in I}p_i$ is the underlying polynomial of a comonoid, then in the corresponding category, $I$ is the set of objects and $p_i$ is the set of outgoing morphisms $i\to\_$. The identities are given by the counit $p\to y$ and the codomains and compositions are given by the comultiplication $p\to p\circ p$.

\begin{example}[Contractible groupoids, $S\yon^S$]\label{ex.contractible}
Let $S$ be a set; the contractible groupoid on $S$ is the category with objects $S$ and a unique morphism $\fun{s}\to \fun{s}'$ for each $\fun{s,s'}\in S$. It corresponds to the comonoid with carrier $S\yon^S$ and counit $S\yon^S\to\yon$ given by evaluation. In other words it is the comonad $\smset\to\smset$ arising from the exponential adjunction for the set $S$. It is often called the \emph{store comonad} in functional programming.
\end{example}

\begin{remark}
In \url{https://www.youtube.com/watch?v=tW6HYnqn6eI}, Richard Garner explains that for any comonoids $C,D$, the $(C,D)$-bimodules in $\poly$ are precisely the \emph{parametric right adjoints} $D\text{-}\smset\to C\text{-}\smset$ between the copresheaf categories.

In \cite{ahman2013distributive} it was shown that the comonoid arising from a distributive law $D\circ C\to C\circ D$ between comonoids in $\poly$ recovers the Zappa-Sz\'ep product \cite{zappa1940sulla} of monoids when the $C,D$ are themselves monoids. 

The point is that comonoids in $\poly$ unexpectedly recover many important notions.
\end{remark}

\chapter{Polynomials and mode-dependent dynamics}\label{chap.poly_mode_dep_dyn}

We next discuss how structures available in $\poly$ describe phenomena in dynamical systems.

\section{Dynamical systems in $\poly$}\label{sec.dynamical}

By a fixed-interface $(A,B)$-dynamical system, we mean a \emph{Moore machine}, i.e.\ a function $\fun{r}\colon S\to B$ (called \emph{readout}), and a function $\fun{u}\colon A\times S\to S$ (called \emph{update}). Given an initial state $\fun{s}_0\in S$, a Moore machine lets us transform any stream $(\fun{a_0},\fun{a_1},\ldots)$ of $A$'s into a stream of $B$'s by repeatedly updating the state:
\[\fun{s}_{n+1}=\fun{u}(\fun{a}_n,\fun{s}_n),\qquad \fun{b}_{n}=\fun{r}(\fun{s}_n).\]

\begin{proposition}\label{prop.tfae}
Let $S,A,B$ be sets. The following are equivalent:
\begin{enumerate}[itemsep=0pt]
	\item Moore machines with inputs $A$, outputs $B$,
	\item coalgebras for the polynomial functor $B\yon^A$,
	\item morphisms in $\poly$ of the form $S\yon^S\to B\yon^A$.
\end{enumerate}
\end{proposition}

The second and third perspectives easily generalize to replacing $B\yon^A$ with an arbitrary polynomial. We prefer the third because it allows us to remain within the category $\poly$, which has such abundant structure. Recall from Example~\ref{ex.contractible} that $S\yon^S$ can be given the structure of a comonoid in $(\poly,\circ,\yon)$, corresponding under Theorem~\ref{thm.cats_comons} to the contractible groupoid on $S$.

\begin{definition}\label{def.mdds}
A \emph{mode-dependent dynamical system} consists of a comonoid $(C,\epsilon,\delta)$ in $(\poly,\circ)$ together with a morphism $f\colon C\to p$ for some polynomial $p$. Here $C$ is called the \emph{state system}, $p$ is called the \emph{interface}, and $f$ is called the \emph{dynamics}.
\end{definition}

We explain the basic idea of Definition~\ref{def.mdds} using the simple comonoid $C=S\yon^S$, and then we'll explain what exactly the comonoid is doing for us. If $p=\sum_{\fun{i}\in p(1)}\yon^{p_\fun{i}}$ is the interface, a morphism $(f,f^\sharp)\colon S\yon^S\to p$ does the following. For each state $\fun{s}$, it returns the current position $f(\fun{s})\in I$; in the case of a monomial $p=B\yon^A$, this $f(\fun{s})\in B$ is the readout for state $\fun{s}$. Then $f^\sharp_{\fun{s}}\colon p_\fun{s}\to S$ sends every direction $\fun{d}\in p_\fun{s}$ to a new state $f^\sharp_{\fun{s}}(\fun{d})$. Again in the case $p=B\yon^A$, the $f^\sharp$ would constitute the update function.

In the morphism $f\colon C\to p$ from Definition~\ref{def.mdds}, the comonoid structure on $C$ provides a canonical morphism $\delta^{n-1}\colon C\to C\circpow{n}$ for each $n$, where $\delta^{-1}=\epsilon$ and $\delta^0=\id$. Since $\circ$ is monoidal, we also have a map $f\circpow{n}\colon C\circpow{n}\to p\circpow{n}$, and composing we obtain
\[C\to C\circpow{n}\to p\circpow{n}.\ 
\footnote{For those readers who are more accustomed to coalgebras, note that one can take the limit of the on-positions functions $C(1)\to p\circpow{n}(1)$, as $n$ increases; this induces the usual map from $C(1)$ to the terminal coalgebra of $p$. In $\poly$ one represents this by a right adjoint $\poly\to\comon(\poly)$ to the forgetful functor. That is, $C\to p$ induces a comonoid morphism $C\nrightarrow\Cat{Cof}(p)$ to the cofree comonoid on $p$, which itself is given by the limit $1\from \yon\mdot p(1)\from \yon\mdot p(\yon\mdot p(1))\from\cdots$ in $\poly$; its set of positions $\Cat{Cof}(p)(1)$ is again the terminal coalgebra on $p$.}
\]
Thus each $\fun{i}\in C(1)$ is endowed with an element of $p\circpow{n}(1)$, which by Eq.~\ref{eqn.strategy} can be understood as a length-$n$ \emph{strategy}
\[
p\circpow n(1)\cong\sum_{\fun{i_1}\in p(1)}\;\prod_{\fun{d_1}\in\fun{p_{i_1}}}\;\sum_{\fun{i_2}\in p(1)}\;\prod_{\fun{d_2}\in\fun{p_{i_2}}}\cdots\;\sum_{\fun{i_n}\in p(1)}\;\prod_{\fun{d_n}\in\fun{p_{i_n}}}1.
\]
It is a choice of a position (move by `player') in $\fun{i}_1\in p(1)$, and for every direction there (move by `opponent') $\fun{d}_1\in\fun{p_{i_1}}$, a choice of position $\fun{i}_2\in p(1)$, etc. Thus a sort of \emph{game} is inherent in the dynamical system itself; it would be interesting to explore a relationship between this and open economic game theory \cite{ghani2016compositional}.

But the map $C\to p\circpow{n}$ does not only give a mapping on positions; it says that for every $n$ choices of directions---each dependent on the last---in $p$, there is a choice of direction, i.e.\ morphism, in the comonoid/category $C$. Thus the history of play is encoded as a morphism in $C$. In the case of coalgebras, where $C=S\yon^S$ is simply a contractible groupoid, there is no information encoded in this history of play, except for its final destination.

\section{Products of interfaces}

The product of polynomials allows one to overlay two different interfaces on the same state system. That is, given dynamical systems $C\to p$ and $C\to q$, there is a unique dynamical system $C\to pq$. This is quite useful for dynamical systems, as we now show.

\begin{example}
Consider two four-state dynamical systems $4\yon^4\to \rr\yon^{\{r,b\}}$ and $4\yon^4\to \rr\yon^{\{g\}}$, each of which gives outputs in $\rr$; we think of $r,b,g$ as red, blue, and green, respectively. We can draw such morphisms as labeled transition systems, e.g.
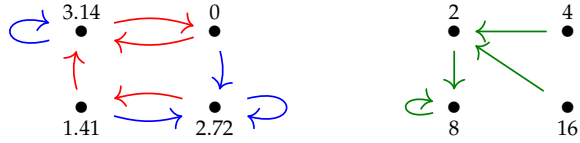
\[
\begin{tikzcd}[row sep=15pt]
	\lmo{3.14}\ar[r, bend left=15pt, red]\ar[loop left=15pt, blue]&
	\lmo{0}\ar[l, bend left=15pt, red]\ar[d, bend left=15pt, blue]\\
	\lmo[under]{1.41}\ar[u,bend left=15pt, red]\ar[r, bend right=15pt, blue]&
	\lmo[under]{2.72}\ar[l, bend right=15pt, red]\ar[loop right=15pt, blue]
\end{tikzcd}
\hspace{.5in}
\begin{tikzcd}[row sep=15pt]
	\lmo{2}\ar[d, green!50!black]&
	\lmo{4}\ar[l, green!50!black]\\
	\lmo[under]{8}\ar[loop left, green!50!black]&
	\lmo[under]{16}\ar[ul, green!50!black]
\end{tikzcd}
\]
Each bullet refers to a state, is labeled by its output position in $\rr$, and has a unique emanating arrow for each sort of input (red and blue, or green), indicating how that state is updated upon encountering said input.

The universal property of products provides a unique way to put these systems together to obtain a morphism $4\yon^4\to(\rr\yon^{\{r,b\}}\times \rr\yon^{\{g\}})=(\rr^2)\yon^{\{r,b,g\}}$. With the examples above, it looks like this:
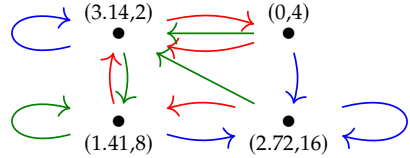
\[
\begin{tikzcd}
	\lmo{(3.14,2)}\ar[r, bend left=15pt, red]\ar[loop left, blue]\ar[d, bend left=15pt, green!50!black]&
	\lmo{(0,4)}\ar[l, bend left=15pt, red]\ar[d, bend left=15pt, blue]\ar[l, green!50!black]\\
	\lmo[under]{(1.41,8)}\ar[u,bend left=15pt, red]\ar[r, bend right=15pt, blue]\ar[loop left, green!50!black]&
	\lmo[under]{(2.72,16)}\ar[l, bend right=15pt, red]\ar[loop right=15pt, blue]\ar[ul, green!50!black]
\end{tikzcd}
\]
Thus the intuitively obvious act of overlaying these dynamical systems falls out of the mathematics, in particular the universal property of products $\times$ in $\poly$. This works for non-monomial (context-dependent) interfaces as well.
\end{example}

\section{Wiring diagrams and mode-dependence}\label{sec.wd_mode_dep}

The Dirichlet product \eqref{eqn.dirichlet} of polynomials and comonoids allows us to juxtapose dynamical systems in an environment. That is, given dynamical systems $C_1\to p_1$ and $C_2\to p_2$, we can form a new dynamical system $(C_1\otimes C_2)\to (p_1\otimes p_2)$.

\begin{example}[Wiring diagrams]\label{ex.wiring_diagram}
Suppose given a wiring diagram such as that in \eqref{eqn.control_diag}; as mentioned in \eqref{eqn.basic_diagram}, the interfaces of the controller and plant are the polynomials $B\yon^C$ and $C\yon^{AB}$, and that of the total system is $C\yon^A$. (Here $A,B,C$ are sets; we will not discuss comonoids again in the remainder of this paper.) All of these are monomials, meaning that the set of directions does not depend on that of positions; this allows us to think of positions as outputs and directions as inputs, drawn on the right and left of boxes respectively. The wiring diagram \eqref{eqn.control_diag} itself is syntax for a morphism
\begin{equation}\label{eqn.wd_control_math}
	B\yon^C\otimes C\yon^{AB}\to C\yon^A.
\end{equation}
On positions the required map $BC\to C$ is the projection, and on directions the required map $BCA\to CAB$ is the obvious symmetry. 
\end{example}

\begin{example}[Mode-dependent wiring diagrams]\label{ex.mode_dep}
In \eqref{eqn.supplier} we depicted a company $C$ changing its supplier of widgets $W$, based on $C$'s internal state. The company was shown with no output wires, but in fact it has two positions corresponding to choosing supplier 1 or supplier 2. Let's redraw it to emphasize its change of position:
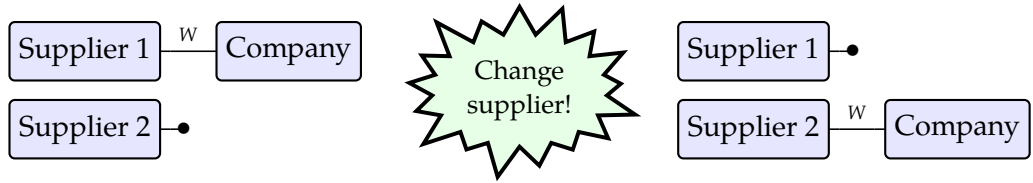
\[
\begin{tikzpicture}[oriented WD, every node/.style={fill=blue!10}]
	\node[bb={0}{1}] (s1) {Supplier 1};
	\node[bb={0}{1}, below=of s1] (s2) {Supplier 2};
	\node[bb={1}{0}, right=0.5 of s1] (c) {Company};
	\draw (s1_out1) to node[above, fill=none, font=\tiny] {$W$} (c_in1);
	\draw (s2_out1) to +(5pt,0) node[fill=none] {$\bullet$};
\begin{scope}[xshift=3.5in]
	\node[bb={0}{1}] (s1') {Supplier 1};
	\node[bb={0}{1}, below=of s1'] (s2') {Supplier 2};
	\node[bb={1}{0}, right=0.5 of s2'] (c') {Company};
	\draw (s2'_out1) to node[above, fill=none, font=\tiny] {$W$} (c'_in1);
	\draw (s1'_out1) to +(5pt,0) node[fill=none] {$\bullet$};
\end{scope}
	\node[starburst, draw, minimum width=2cm, minimum height=2cm,align=center,fill=green!10, font=\small,line width=1.5pt] at ($(c.east)!.5!(s2'.west)$)
{Change\\supplier!};
\end{tikzpicture}
\]
The company has interface $2\yon^W$, and the each supplier has interface $W\yon$; let's take the total system interface (undrawn) to be the closed system $\yon$. Then this mode-dependent wiring diagram is just a map $2\yon^W\otimes W\yon\otimes W\yon\to\yon$. Its on-positions function $2W^2\to1$ is uniquely determined, and its on-directions function $2W^2\to W$ is the evaluation. In other words, the company's position determines which supplier from which it receives widgets.

Similarly we could say that the person in \eqref{eqn.someone} has interface $2\yon$, the units have interfaces $X\yon$ and $\yon^X$ respectively, and the whole system is closed; that is, the diagram represents a morphism $2\yon\otimes X\yon\otimes \yon^X\to\yon$. We did not mention but need that unit B has a default value, say $\fun{x_0}\in X$, for when its input wire is unattached. The morphism $2X\yon^X\to\yon$ is uniquely determined on positions, and on directions it is given by cases $(1,x)\mapsto x_0$ and $(2,x)\mapsto x$.
\end{example}

\printbibliography

\end{document}